\documentclass{amsproc}
\newtheorem{theorem}{Theorem}[section]
\newtheorem{lemma}[theorem]{Lemma}

\theoremstyle{definition}
\newtheorem{definition}[theorem]{Definition}

\theoremstyle{remark}

\numberwithin{equation}{section}

\usepackage{amssymb}
\usepackage[latin1]{inputenc} %Neue Codepage
\usepackage{textcomp}

\newcommand{\nz}{\mathbb{N}} % Nat"urliche Zahlen
\newcommand{\N}{\nz}
\newcommand{\rz}{\mathbb{R}} % Relle Zahlen
\newcommand{\grad}{\nabla}

\newcommand{\di}{\,{\rm d}}

\numberwithin{equation}{section}

\makeindex
\usepackage{pstricks}

\newcommand{\pstheight}{0.8}

\newcommand{\pstpic}[1]{%
\pspicture(0.000000,-0.130000)(0.9470,\pstheight)%
\ifx\nofigs\undefined%
\catcode`@=11%
\rput[lt](\pstxo,\pstyo){%
\begin{minipage}{3cm}%
\footnotesize%
#1
\end{minipage}}%
}
\newcommand{\pstpend}{%
\catcode`@=12%
\fi%
\endpspicture%
}
\newcommand{\numberofit}{\text{iter\#}}

\newcommand{\bil}{\Gamma}
\newcommand{\abil}{\tilde{\Gamma}}
\newcommand{\R}{\mathbb{R}}
\newcommand{\q}{\boldsymbol{q}}
\newcommand{\Q}{\boldsymbol{Q}}
\begin{document}

% \title[short text for running head]{full title}
\title{Recovering boundary conditions in inverse Sturm-Liouville problems}

%    Only \author and \address are required; other information is
%    optional.  Remove any unused author tags.

%    author one information
% \author[short version for running head]{name for top of paper}
\author{Norbert Röhrl}
\address{Fachbereich Mathematik, Univerität Stuttgart, 70550 Stuttgart, Germany}
\curraddr{}
\email{norbert.roehrl@mathematik.uni-stuttgart.de}
\thanks{}

%    \subjclass is required.

\date{}

\dedicatory{}
%    "Communicated by" -- provide editor's name; required.
%\commby{}

% Abstract is required.  

\begin{abstract} We introduce a variational algorithm, which solves
  the classical inverse Sturm-Liouville problem when two spectra are
  given. In contrast to other approaches, it recovers the potential as
  well as the boundary conditions without a priori knowledge of the
  mean of the potential.  Numerical examples show that the algorithm
  works quite reliable, even in the presence of noise. A proof of the
  absence of strict local minimizers of the functional supports the
  observation, that a good initial guess is not essential.
\end{abstract}

\maketitle

%    Text of article.
\section{Introduction}
The inverse Sturm-Liouville problem was first systematically studied
by Borg in 1946 \cite{Borg1946}. He already proved that all
information needed to reconstruct the potential is in  two sequences of
eigenvalues, and applied this to the question if one could hear the
mass density of a guitar string.

With modern computers, the question for efficient algorithms to
actually compute the potential gained importance
\cite{McLaughlin1986}. This work was inspired by two different
approaches to this problem.

One was developed by Rundell and Sacks \cite{RundellSacks1992}
and is based on the Gelfand-Levitan-Marchenko kernels. It
is elegant and efficient, but also invariably needs the mean $\int_0^1Q
\di x$ of the potential and the boundary conditions as additional
inputs besides the two spectra.  The second method is variational and
was created by Brown, Samko, Knowles, and Marletta
\cite{BrownSamkoetal2003}. It does not need the mean as a separate
input, but it is unknown if it also can be used to recover the
boundary conditions.

We also want to mention a related recovery method by Lowe, Pilant, and
Rundell \cite{Loweetal1992}, which uses a finite basis ansatz, and
solves the inverse problem by Newton's method without requiring the
mean as input.

In this paper we extend the variational method we introduced in
\cite{Roehrl2005} to recover potential and boundary conditions in the
case when only two finite sequences of eigenvalues
are given.  Having less reliable information it is not as robust under
noisy input, but we still get reasonable results.

In the following section we will define the functional and
exhibit some essential properties. The numerical examples are
discussed in the third section and section four finally contains
the proof of the absence of strict local minimizers.

\section{Definition and Properties of the Functional}\label{sec:def}
\noindent We consider the Sturm-Liouville equation
\renewcommand{\theequation}{SL}
\begin{equation}
\label{SL}
{-u'' + q(x)u = \lambda u}
\end{equation}
on $[0,1]$ with $q(x)\in L^2([0,1],\rz)$ real, and separated boundary
conditions
\renewcommand{\theequation}{$h_0h_1$}
\begin{equation}
%\tag{$h_0h_1$}
 h_0 u(0) +u'(0)=0,\quad h_1u(1) +u'(1) = 0\,.
\end{equation}
\renewcommand{\theequation}{\arabic{section}.\arabic{equation}}

The corresponding eigenvalues satisfy the asymptotic formula
\cite{IsaacsonTrubowitz1983}
\begin{equation}\label{eigenasymptotics}
 \lambda_{n}=\pi^2n^2+ 2(h_1-h_0) + \int_0^1 q(s) \di s + a_n\,, 
\end{equation}
where $(a_n)\in l^2$.  It is a classical result \cite{Borg1946,Levinson1949},
that the potential $q$ is uniquely determined by two sequences of
eigenvalues corresponding to boundary conditions $(h_0h_1)$ and
$(h_0h_2)$ with $h_1\neq h_2$. Moreover it can be shown, that those
sequences also uniquely determine the boundary conditions
\cite{Levitan1987}.

Therefore, two sequences of eigenvalues contain all information
necessary to recover the potential as well as the corresponding
boundary conditions. For notational convenience we write the parameters
of the two Sturm-Liouville problems as vectors
\[
\boldsymbol{q_1}:=(h_0,h_1,q)\,, \quad \boldsymbol{q_2}:=(h_0,h_2,q)\,,
\]
$\lambda_{\boldsymbol{q}_i,n}$ for the $n$-th eigenvalue of problem $\q_i$, and 
\[
\q := (h_0,h_1,h_2,q)
\]
for the full problem.

Now we define a least squares functional on the eigenvalues, which has
the solution of the inverse problem as zero.  
\begin{definition}
  Suppose we are given (partial) spectral data $\lambda_{\Q_i,n}$ with
  $(i,n)$ in $I \subseteq \{1,2\}× \N$ of an unknown Sturm-Liouville problem
  $\Q=(H_0,H_1,H_2,Q)$.  For a trial problem $\q$ and positive weights
  $\omega_{i,n}$, we define the functional
\begin{equation}\label{functional}
G(\boldsymbol{q}) := \sum_{(i,n)\in I} \omega_{i,n}(\lambda_{\boldsymbol{q}_i,n} - \lambda_{\Q_i,n})^2 \,.
\end{equation}
\end{definition}

We note that $G(\q)$ is positive, and zero if and only if both given
sequences of eigenvalues match those of $\q$. If we have full
knowledge of the two sequences $\lambda_{\Q_i,n}, (i,n)\in \{1,2\}× \N$, this
determines $\q$ uniquely, and hence $\q=\Q$.

To find such a $\q$, we minimize the functional with a conjugate
gradient descent algorithm. First, for numerical stability it is good
to know, that for each pair of interlacing sequences, 
\[
   \lambda_{1,n}<\lambda_{2,n}<\lambda_{1,n+1}\qquad \mathrm{ or } \qquad \lambda_{2,n}<\lambda_{1,n}<\lambda_{2,n+1}\,,
\]
which satisfy
the asymptotics \eqref{eigenasymptotics}, there is a $\q$, with
$G(\q)=0$ \cite{Levitan1987}.
 
The gradient of $\lambda_{\boldsymbol{q}_i,n}$ wrt.~$\boldsymbol{q}_i$ is \cite{IsaacsonTrubowitz1983}
\[
\grad {\lambda}_{\boldsymbol{q}_i,n}= 
\begin{pmatrix}
\frac{\partial{\lambda}_{\boldsymbol{q}_i,n}}{\partial h_0}\\
\frac{\partial{\lambda}_{\boldsymbol{q}_i,n}}{\partial h_1}\\
\frac{\partial{\lambda}_{\boldsymbol{q}_i,n}}{\partial h_2}\\
\frac{\partial{\lambda}_{\boldsymbol{q}_i,n}}{\partial q}\\
\end{pmatrix}=
\begin{pmatrix}
  -g^2_{\boldsymbol{q}_i,n}(0)\\
  g^2_{\boldsymbol{q}_i,n} (1) \delta_{i,1}\\
  g^2_{\boldsymbol{q}_i,n} (1) \delta_{i,2}\\
  g^2_{\boldsymbol{q}_i,n} (x)
\end{pmatrix}\,,
\]
where $g_{\q_i,n}$ denotes the eigenfunction
corresponding to $\lambda_{\q_i,n}$ with $\|g_{\q_i,n}\|=1$.

It follows that the gradient of the functional is given by
\[
\grad{G}(\q) = 2  \sum_{(i,n)\in I} \omega_{i,n}(\lambda_{\boldsymbol{q}_i,n} - \lambda_{\Q_i,n})\grad {\lambda}_{\boldsymbol{q}_i,n}\,,
\]
if  $(\omega_{i,n})$ is summable. 

Theorem~\ref{linind} below shows that the gradients $\grad \lambda_{\q_i,n}$
are linearly independent in $\R^3× L^2(0,1)$. This immediately implies the
essential convexity of the functional:

\begin{theorem}\label{mainth}
  Let $\boldsymbol{q}_1=(h_0,h_1,q)$ and $\boldsymbol{q}_2=(h_0,h_2,q)$
  be two Sturm-Liouville problems with $h_1\neq h_2$. If $I$ is finite or
  $(\omega_{i,n})$ is summable, the functional $G$ has no local minima at
  $\q$ with $G(\q)>0$, i.e.
\[
\grad {G}(\q) = 0  \Longleftrightarrow G(\q)=0\,.
\]
\end{theorem}
Thus a conjugate gradient algorithm will not get trapped in local
minima, as we will also observe in the examples. 

\section{Numerical Examples}\label{sec:examples}
We use the standard Polak-Ribiere conjugate gradient descent algorithm
\cite{NumericalRecipes} to approximate the gradient flow and thus
minimize the functional. To give the basic idea, we explain the
simpler steepest descent:
\begin{enumerate}
\item choose initial potential and boundary conditions $\q^{(0)}=(h_0^{(0)},h_1^{(0)},h_2^{(0)},q^{(0)})$
\item while $G(\q^{(j)})$ too big do
  \begin{enumerate}
  \item compute the gradient $\grad G(\q^{(j)})$
  \item minimize the one dimensional function $G(\q^{(j)}) - \alpha \grad G(\q^{(j)})$ wrt.~$\alpha$  
  \item set $\q^{(j+1)}$ equal to the minimizing potential 
  \end{enumerate}
\end{enumerate}
This straight forward minimization scheme has a major disadvantage:
consecutive gradients are always orthogonal. To avoid this, conjugate
gradient descent computes the direction for the one dimensional
minimization using the current and previous gradients. 

Note that the boundary points $g_{\q_i,n}(0)$ and $g_{\q_i,n}(1)$, needed in
the computation of the gradient $\grad \lambda_{\q_i,n}$, can be computed in
a numerically well behaved way. Given the eigenvalue $\lambda_{\q_i,n}$, we
can compute a multiple of the eigenfunction by solving an initial value
problem. The value of the eigenfunction at the boundary then is
just the (exactly known) initial value divided by the $L^2$ norm of the
initial value solution.

We first apply the algorithm to a popular non-continuous potential \cite{BrownSamkoetal2003,Roehrl2005}
\[
Q(x)= (7x-0.7)\chi_{(0.1,0.3]}(x) + (3.5-7x)\chi_{(0.3,0.5]}(x) + 4\chi_{(0.7,0.9]}(x) + 2\chi_{(x,0.9,1]}(x)\,,
\] 
and choose
\begin{equation}\label{sbc}
\Q=(3,3,0,Q) \qquad \text{ and } \qquad \q^{(0)}=(2,4,-1,0)
\end{equation}
with $I=\{1,2\}× \{0,\dots,29\}$ and weights $\omega_{i,n}=1$. This will be our
default setting, unless noted otherwise. In the figures we write
$\Delta_2=\|q-Q\|_2$ for the $L^2$ error of the reconstruction and
$bc=(h_0;h_1;h_2)$ for the boundary conditions of the current
approximation.

\input{gnuplothead.tex}

\newcommand{\pstxo}{-3.2}
\newcommand{\pstyo}{4.8}
\newcommand{\Mynumb}{%
\[\arraycolsep0.1cm\begin{array}{rcl}
    G(\q) &=& 8.6764\cdot10^{-5}\\
    \Delta_2 &=& 0.66227\\
    \numberofit &=& 150\\
    \text{bc} &=& (2.844; 3.035; 0.033)
  \end{array}\]}
\newcommand{\Mynumbr}{}
\newcommand{\Mynumbb}{%
$G(\q)$ over iter\# 
}
\newcommand{\Mynumbbr}{%
$\Delta_2$ over iter\# 
}

\begin{figure}
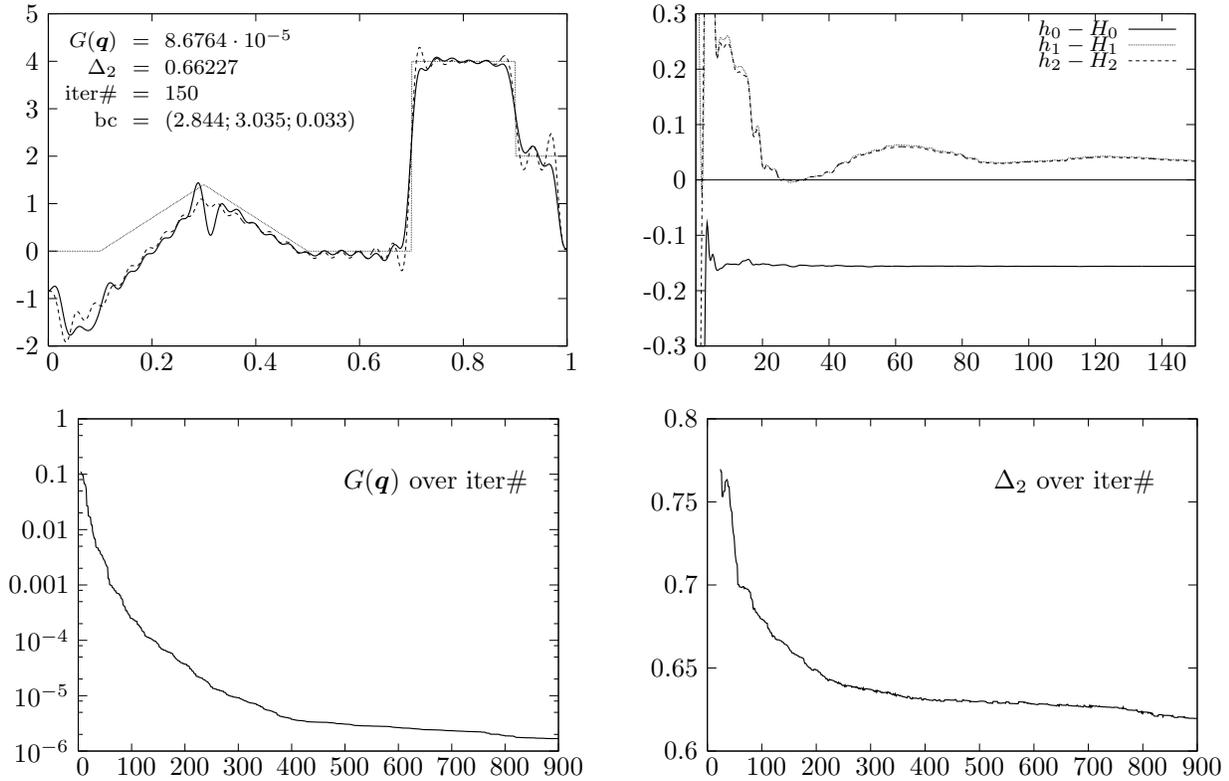
%
\begin{tabular}{cc}%
\unitlength1cm%
\begin{picture}(8,5)%
\pstpic{\Mynumb}%
\psset{xunit=8.0cm,yunit=5.0cm}%
\rput(-.5,0){\input{norm_1.tex}}%
\pstpend%
\end{picture}%
&%
\begin{picture}(8,5)%
\pstpic{\Mynumbr}%
\psset{xunit=8.0cm,yunit=5.0cm}%
\rput(-.5,0){\input{norm_1_abc.tex}}%
\pstpend%
\end{picture}\\[4.5cm]%
\begin{picture}(8,5)%
%\psgrid(-10,-10)(10,10)
\rput(-7.5,0){%
\pstpic{}%\Mynumbb}%
\rput(6,4){\Mynumbb}%
\psset{xunit=8.0cm,yunit=5.0cm}%
\input{norm_1_convG.tex}\pstpend%
}%
\end{picture}%
&%
\begin{picture}(8,5)%
\rput(-3.5,0){%
\pstpic{}%\Mynumbbr}%
\rput(6,4){\Mynumbbr}%
\psset{xunit=8.0cm,yunit=5.0cm}%
\input{norm_1_convL2.tex}\pstpend%
}%
\end{picture}%
\end{tabular}%
\caption{Graph of $\q^{(150)}$,
  $\q^{(840)}$ (light) and boundary
  conditions, $G(\q)$, and $L^2$ error versus the number of iterations.}
\label{fig1}
\medskip
\end{figure}%

% frprmn:8.67636568297e-05
% P00150
% D00150
% DP01939
% ERROR: l2norm not within tolerance!
% l2error 0.662273178786
% diff in mean0.36570885
% poten_real( 2.84410797236; 3.03500671043; 0.0334927100709)

% frprmn:1.7568305003e-06
% P00100
% D00100 iteration 840
% DP01513
% l2error 0.62091127375
% diff in mean0.3457922019
% poten_real( 2.8469410386; 3.02742210721; 0.0263790195669)

The results (figure~\ref{fig1}) are comparable to the
alternative boundary condition example with given boundary conditions
in \cite{Roehrl2005}. We get quite good results at 150 iterations, which
keep getting better as we minimize the functional.

From around 150 iterations the boundary conditions stay almost
constant. This suggests to reset the trial potential $q^{(j)}$ to zero
at some iteration number $j$, while keeping the boundary conditions.
In other words, we reset the forth component of
$\q^{(j)}=(h_0^{(j)},h_1^{(j)},h_2^{(j)},q^{(j)})$, and keep the
others fixed. In figure \ref{fig3} we indeed attain a significantly
better approximation by setting the potential to zero a couple of
times. In all our examples this was a very useful strategy to get
faster convergence. But it is just heuristics -- we do not really know
how to choose the optimal number of iterations $j$. In practice, we
wait until the boundary conditions stabilize and then set the
potential to zero. This can be repeated until the convergence speed of
the functional $G$ does not improve any more.

The graph in figure \ref{fig3b} demonstrates, that the
reconstruction of a smooth potential, using the same boundary
conditions and number of eigenvalues, yields more accurate results. 

\renewcommand{\pstxo}{0.8}
\renewcommand{\pstyo}{4.8}
\renewcommand{\Mynumb}{%
\[\arraycolsep0.1cm\begin{array}{rcl}
    G(\q) &=& 3.62063\cdot10^{-6}\\
    \Delta_2 &=& 0.404896\\
    \numberofit &=& 350\\
    \text{bc} &=& (2.9322; 3.0285; 0.02742)
  \end{array}\]}

\begin{figure}%
\unitlength1cm%
\begin{center}%
\begin{picture}(8,5)%
\pstpic{\Mynumb}%
\psset{xunit=8.0cm,yunit=5.0cm}%
\rput(0,0){\input{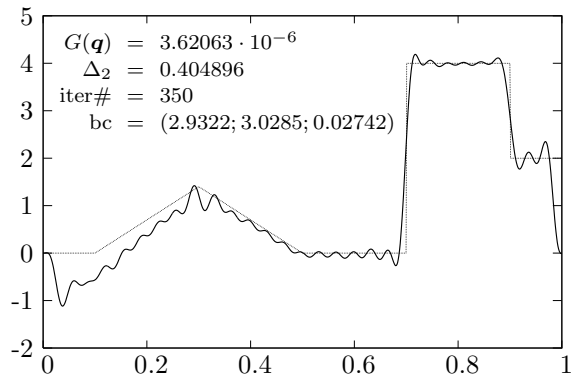}}%
\pstpend%
\end{picture}%
\end{center}%
%\begin{caption}
\vspace{-.5cm}%
\caption{Above example with setting
  $q^{(30)}=q^{(60)}=q^{(110)}=0$.}%
\label{fig3}%
%\end{caption}
\end{figure}%

\renewcommand{\Mynumb}{
\[\ \ \arraycolsep0.1cm\begin{array}{rcl}
    G(\q) &=& 1.3369\cdot10^{-8}\\
    \Delta_2 &=& 0.0387494\\
    \numberofit &=& 1887\\
    \text{bc} &=& (2.99; 3.00;\\ 
    & & 3.99\cdot 10^{-4})
  \end{array}\]}
\renewcommand{\pstxo}{0.9}

\begin{figure}%
\unitlength1cm%
\begin{center}%
\begin{picture}(8,5)%
\pstpic{\Mynumb}%
\psset{xunit=8.0cm,yunit=5.0cm}%
\rput(0,0){\input{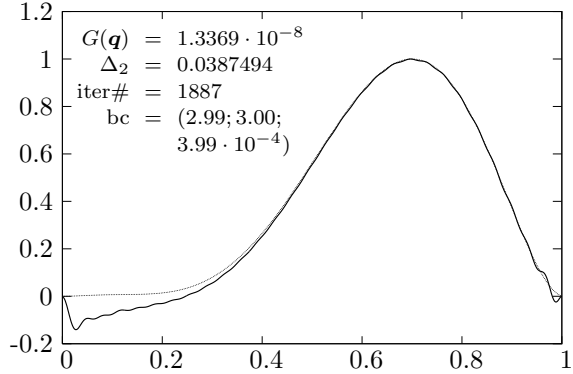}}%
\pstpend%
\end{picture}%
\end{center}%
%\begin{caption}
\vspace{-.5cm}%
\caption{Approximation of a smooth
  potential.}%
\label{fig3b}%
%\end{caption}
\end{figure}%

% frprmn:1.33686563992e-08| cont 822| cont 101| => total 1887
% P00964^@
% D00964^@
% DP17783^@
% l2error 0.0387494363553
% diff in mean0.0226910310875
% poten_real( 2.98950563116; 3.00042636159; 0.000399501479136)
% Newode 0x8191330 0x8191358

% frprmn:3.62063003676e-06
% P00240^@
% D00240^@
% DP03216^@
% l2error 0.4048955394
% diff in mean0.17548619371
% poten_real( 2.9321947582; 3.02847860199; 0.0274170902343)

This overall behavior is also true for worse guesses of the initial
boundary conditions. If we take for example 
\[
\q^{(0)}=(0,0,0,0)\,,
\]
the boundary conditions converge slowly, but steadily (figure~\ref{fig4}).
Again, setting $q^{(j)}=0$ a couple of times increases the speed of convergence 
dramatically.

\renewcommand{\pstxo}{0}
\renewcommand{\pstyo}{4.0}%4.8
\renewcommand{\Mynumb}{%
}
\renewcommand{\Mynumbr}{%
\normalsize $G(\q)$ over iter\# 
}
\begin{figure}
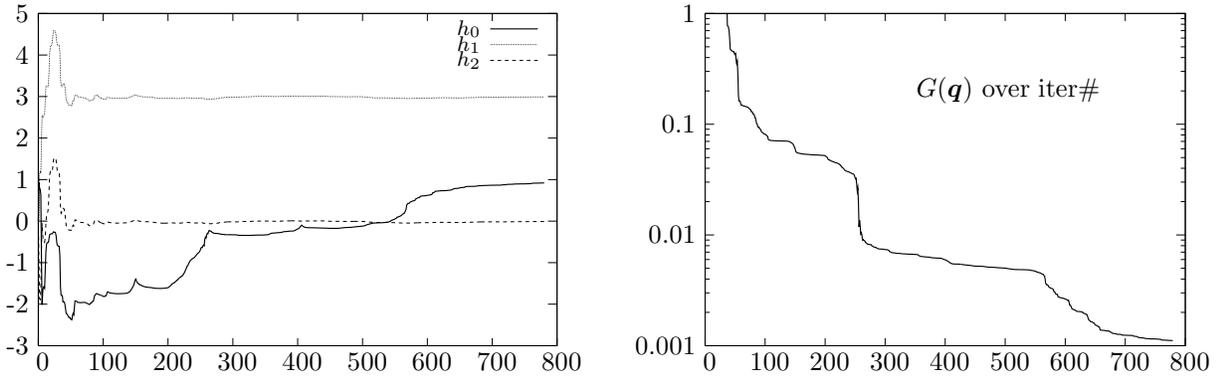
%
\begin{tabular}{cc}%
\unitlength1cm%
\begin{picture}(8,5)%
\pstpic{\Mynumb}%
\psset{xunit=8.0cm,yunit=5.0cm}%
\rput(-.5,0){\input{000_abc.tex}}%
\pstpend%
\end{picture}%
&%
\begin{picture}(8,5)%
\pstpic{\Mynumbr}%
\psset{xunit=8.0cm,yunit=5.0cm}%
\rput(-.5,0){\input{000_convG.tex}}%
\pstpend%
\end{picture}%
\end{tabular}%
%\begin{caption}
\vspace{-.5cm}
\caption{The convergence of the boundary
conditions and of the functional for initial problem $\q^{(0)}=(0,0,0,0)$.}
\label{fig4}
%\end{caption}
\end{figure}%

% frprmn:0.0217137582634
% P00023^@
% D00023^@
% DP00290^@
% l2error 0.567095400432
% diff in mean0.217410247185
% poten_real( 2.92777641548; 3.06659313072; 0.0661009931606)

% frprmn:0.0160616321683
% P00055^@
% D00055^@
% DP00692^@
% ERROR: l2norm not within tolerance!
% l2error 1.16280623807
% diff in mean0.238631449693
% poten_real( 2.92155426183; 3.07069601501; 0.0686959605872)

In the case of noisy data, the algorithm is of course much more
unstable than the version with fixed boundary conditions
\cite{Roehrl2005}. For the following examples, we add random noise
$|\tilde{\lambda}_{\Q_i,n} - \lambda_{\Q_i,n}|\leq r$ to the eigenvalues. To see the
limitations of this approach, we first use $r=0.1$ and
$H_0=h_0=H_1=h_1=3$, $H_2=h_2=0$ (figure~\ref{fig5}), i.e.~we already
start with the correct boundary conditions.

\renewcommand{\pstxo}{-3.2}
\renewcommand{\pstyo}{4.8}
\renewcommand{\Mynumb}{%
\[\arraycolsep0.1cm\begin{array}{rcl}
    G(\q) &=& :0.021714\\
    \Delta_2 &=& 0.567095\\
    \numberofit &=& 23\\
    \text{bc} &=& (2.928; 3.067; 0.0661)
  \end{array}\]}
\renewcommand{\Mynumbr}{}

\begin{figure}
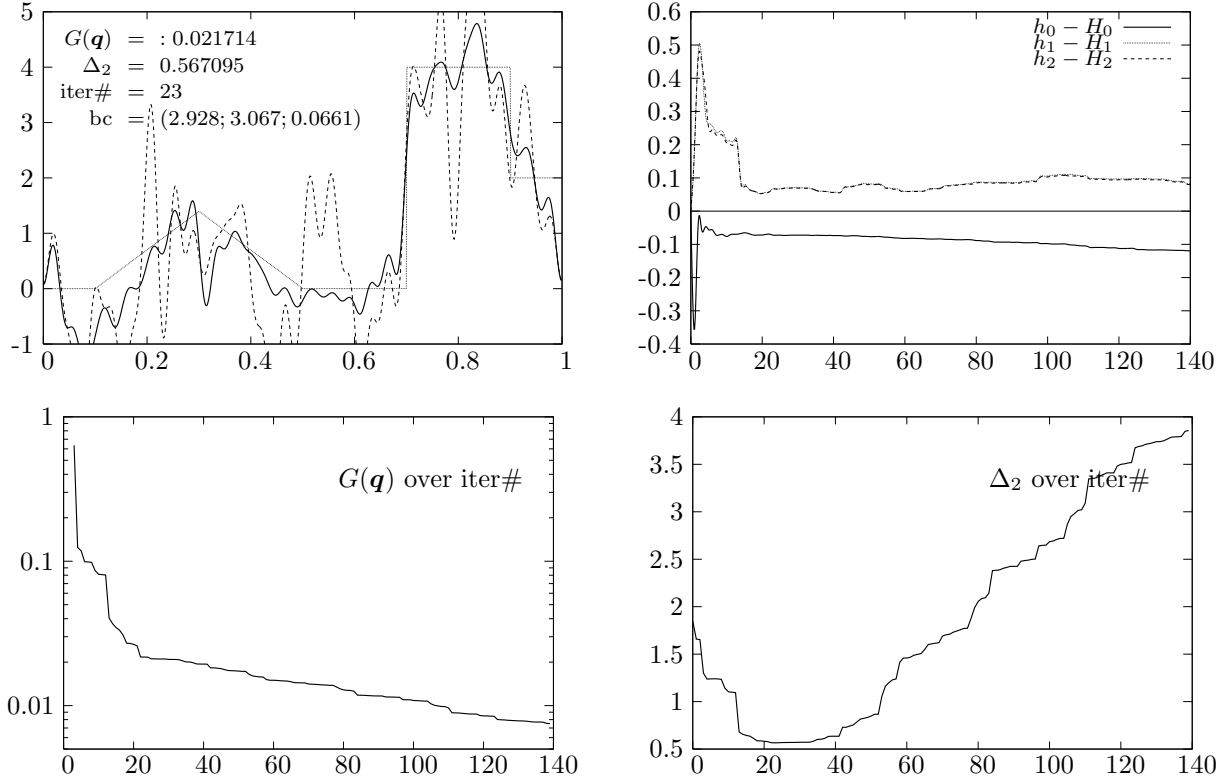
%
\begin{tabular}{cc}%
\unitlength1cm%
\begin{picture}(8,5)%
\pstpic{\Mynumb}%
\psset{xunit=8.0cm,yunit=5.0cm}%
\rput(-.5,0){\input{r0.1_330_330.tex}}%
\pstpend%
\end{picture}%
&%
\begin{picture}(8,5)%
\pstpic{\Mynumbr}%
\psset{xunit=8.0cm,yunit=5.0cm}%
\rput(-.5,0){\input{r0.1_330_330_abc.tex}}%
\pstpend%
\end{picture}\\[4.5cm]%
\begin{picture}(8,5)%
%\psgrid(-10,-10)(10,10)
\rput(-7.5,0){%
\pstpic{}%\Mynumbb}%
\rput(6,4){\Mynumbb}%
\psset{xunit=8.0cm,yunit=5.0cm}%
\input{r0.1_330_330_convG.tex}\pstpend%
}%
\end{picture}%
&%
\begin{picture}(8,5)%
\rput(-3.5,0){%
\pstpic{}%\Mynumbbr}%
\rput(6,4){\Mynumbbr}%
\psset{xunit=8.0cm,yunit=5.0cm}%
\input{r0.1_330_330_convL2.tex}\pstpend%
}%
\end{picture}%
\end{tabular}%
\caption{Example with
  errors of magnitude $r=0.1$. The best $L^2$ approximation at 23 
  iterations and $q^{(55)}$.}
\label{fig5}
\medskip
\end{figure}%

% \renewcommand{\pstxo}{0} \renewcommand{\pstyo}{4.0}%4.8
% \renewcommand{\Mynumb}{%
%   $G(\q)$ over iter\# } 
% \renewcommand{\Mynumbr}{%
%   $\Delta_2$ over iter\# }
% \pstfig{r0.1_330_330_convG.tex}{r0.1_330_330_convL2.tex}{fig6}{Values
%   of the functional $G(\q)$ resp.~ $L^2$ error $\Delta_2$ versus number of
%   iterations.}

The algorithm passes through a potential, which is reasonably close to
the original potential $Q$. But from there, the steepest descent leads
to a potential, which is not in any way similar to the one we want to
recover. In the graph of the $L^2$ error, we see that there are
roughly 20 iterations of good approximations. Afterwards the
approximations quickly get worse than our initial guess.

Yet,   for smaller  errors in  the   eigenvalues, this effect  is less
dramatic.  Setting $r=0.01$ and  using the problem \eqref{sbc}, we get
good  approximations for   around 100  iterations (figure~\ref{fig7}).
The $L^2$ error also  rises more slowly than  for the case with larger
noise level.

% frprmn:0.000139224899859
% P00089
% D00089
% DP01216
% ERROR: l2norm not within tolerance!
% l2error 0.669477733816
% diff in mean0.36353038757
% poten_real( 2.84323267354; 3.03428926576; 0.0326983301928)

% frprmn:2.0246191731e-05
% P00300^@
% D00300^@
% DP04624^@
% ERROR: l2norm not within tolerance!
% l2error 0.852575415044
% diff in mean0.388343922701
% poten_real( 2.83021135675; 3.03412402575; 0.0329349947799)

\renewcommand{\pstxo}{-3.2}
\renewcommand{\pstyo}{4.8}
\renewcommand{\Mynumb}{%
\[\arraycolsep0.1cm\begin{array}{rcl}
    G(\q) &=& 1.39225\cdot10^{-4}\\
    \Delta_2 &=&  0.669478\\
    \numberofit &=& 89\\
    \text{bc} &=& ( 2.8302; 3.0343; 0.0329)
  \end{array}\]}
\renewcommand{\Mynumbr}{}
\begin{figure}
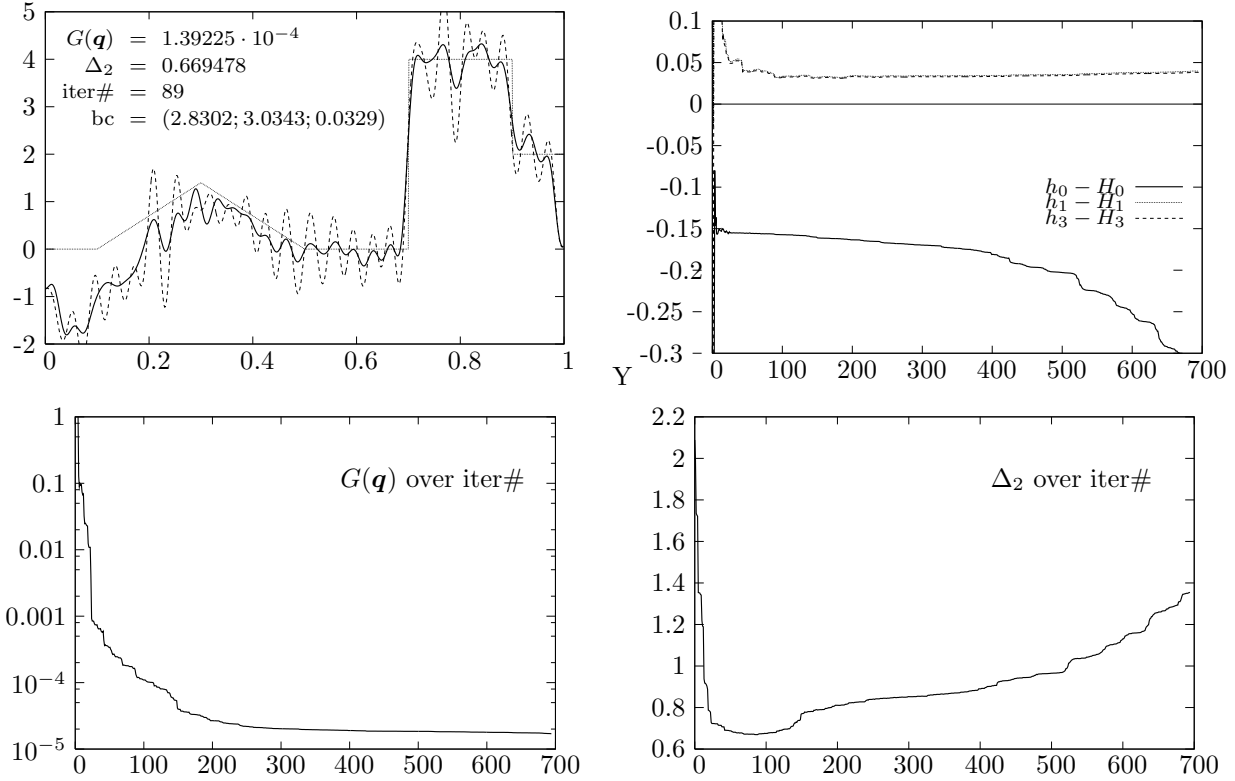
%
\begin{tabular}{cc}%
\unitlength1cm%
\begin{picture}(8,5)%
\pstpic{\Mynumb}%
\psset{xunit=8.0cm,yunit=5.0cm}%
\rput(-.5,0){\input{r0.01_II.tex}}%
\pstpend%
\end{picture}%
&%
\begin{picture}(8,5)%
\pstpic{\Mynumbr}%
\psset{xunit=8.0cm,yunit=5.0cm}%
\rput(-.5,0){\input{r0.01_330_24-1_abc}}%
\pstpend%
\end{picture}\\[4.5cm]%
\begin{picture}(8,5)%
%\psgrid(-10,-10)(10,10)
\rput(-7.5,0){%
\pstpic{}%\Mynumbb}%
\rput(6,4){\Mynumbb}%
\psset{xunit=8.0cm,yunit=5.0cm}%
\input{r0.01_330_24-1_convG.tex}\pstpend%
}%
\end{picture}%
&%
\begin{picture}(8,5)%
\rput(-3.5,0){%
\pstpic{}%\Mynumbbr}%
\rput(6,4){\Mynumbbr}%
\psset{xunit=8.0cm,yunit=5.0cm}%
\input{r0.01_330_24-1_convL2}\pstpend%
}%
\end{picture}%
\end{tabular}%
\caption{Example with
  errors of magnitude $r=0.01$. The best $L^2$ approximation at 89
  iterations and $q^{(300)}$.}
\label{fig7}
\medskip
\end{figure}%

% \renewcommand{\pstxo}{0} \renewcommand{\pstyo}{4.0}%4.8
% \renewcommand{\Mynumb}{%
%   $G(\q)$ over iter\# } 
% \renewcommand{\Mynumbr}{%
%   $\Delta_2$ over iter\# }
% \pstfig{r0.01_330_24-1_convG.tex}{r0.01_330_24-1_convL2}{fig8}{Values
%   of the functional $G(\q)$ resp.~ $L^2$ error $\Delta_2$ versus number of
%   iterations.}

Setting $q^{(89)}=0$ again improves the performance significantly
(figure~\ref{fig9}). We get reasonable approximations for all 400
iterations but the first 20 after each setting the potential to zero.

% frprmn:0.000150128052323
% P00148^@
% D00148^@
% DP01900^@
% l2error 0.471443957363
% diff in mean0.203071210395
% poten_real( 2.92294340583; 3.03446860224; 0.0334690789408)

% frprmn:6.35119066412e-05
% P00300^@
% D00300^@
% DP04013^@
% l2error 0.520395527551
% diff in mean0.211172485123
% poten_real( 2.92023643647; 3.03597070087; 0.0346928363164)

\renewcommand{\pstxo}{-3.2}
\renewcommand{\pstyo}{4.8}
\renewcommand{\Mynumb}{%
\[\arraycolsep0.1cm\begin{array}{rcl}
    G(\q) &=&1.50128 \cdot10^{-4} \\
    \Delta_2 &=&  0.47144\\
    \numberofit &=& 237\\
    \text{bc} &=& ( 2.9229; 3.0345; 0.03347)
  \end{array}\]}
\renewcommand{\Mynumbr}{}
\begin{figure}%
\begin{tabular}{cc}%
\unitlength1cm%
\begin{picture}(8,5)%
\pstpic{\Mynumb}%
\psset{xunit=8.0cm,yunit=5.0cm}%
\rput(-.5,0){\input{r0.01c_II.tex}}%
\pstpend%
\end{picture}%
&%
\begin{picture}(8,5)%
\pstpic{\Mynumbr}%
\psset{xunit=8.0cm,yunit=5.0cm}%
\rput(-.5,0){\input{r0.01c_abc.tex}}%
\pstpend%
\end{picture}\\[4.5cm]%
\begin{picture}(8,5)%
%\psgrid(-10,-10)(10,10)
\rput(-7.5,0){%
\pstpic{}%\Mynumbb}%
\rput(6,4){\Mynumbb}%
\psset{xunit=8.0cm,yunit=5.0cm}%
\input{r0.01c_convG.tex}\pstpend%
}%
\end{picture}%
&%
\begin{picture}(8,5)%
\rput(-3.5,0){%
\pstpic{}%\Mynumbbr}%
\rput(6,4){\Mynumbbr}%
\psset{xunit=8.0cm,yunit=5.0cm}%
\input{r0.01c_convL2.tex}\pstpend%
}%
\end{picture}%
\end{tabular}%
\caption{Example with
  errors of magnitude $r=0.01$ and $\q^{(89)}$ set to zero. The best $L^2$ approximation at 237
  iterations and $q^{(389)}$.}
\label{fig9}
\medskip
\end{figure}%

% \renewcommand{\pstxo}{0} \renewcommand{\pstyo}{4.0}%4.8
% \renewcommand{\Mynumb}{%
%   $G(\q)$ over iter\# } 
% \renewcommand{\Mynumbr}{%
%   $\Delta_2$ over iter\# }
% \pstfig{r0.01c_convG.tex}{r0.01c_convL2.tex}{fig10}{Values
%   of the functional $G(\q)$ resp.~ $L^2$ error $\Delta_2$ versus number of
%   iterations.}

\section{Linear Independence of the Gradients}

First, we borrow the central lemma of the independence proof of the original 
functional \cite{Roehrl2005}. To that end, we define the Wronskian
 $[f,g]= fg'-f'g$ and the bilinear form
\[
\begin{array}{rrcl}
\bil: &H^1([0,1],\R)^2& \longrightarrow & \R\\
   &   (f,g) &\mapsto &  \int_0^1 [f,g] \di x
\end{array}\,,
\]
which is bounded by 
\[
|\bil(f,g)| \leq \|f\|_{H_1} \|g\|_{H_1},\text{ i.e. } \|\bil(f,\cdot) \| = \|f\|_{H_1}\,. 
\]
(We use the definition $\|f\|_{H^1}=\sqrt{\|f\|_{L^2}^2+ \|f'\|_{L^2}^2}$ with
distributional derivatives.) 

Let $s_{i,n,q}$ and $c_{i,n,q}$ be the solutions
of the differential equation \eqref{SL} for the eigenvalue parameter
$\lambda_{\q_i,n}$ and initial values
\[
\begin{array}{rclcrcl}
  s_{i,n,q}(1) &=& 1\,,&\phantom{---}&\qquad c_{i,n,q}(1) &=& 1\,, \\
  s_{i,n,q}'(1) &=& -h_1\,,& &\qquad c_{i,n,q}'(1) &=& -h_2\,.
\end{array}
\]
\begin{lemma}[\cite{Roehrl2005}]
  Given two Sturm-Liouville problems $\boldsymbol{q}_1=(h_0,h_1,q)$
  and $\boldsymbol{q}_2=(h_0,h_2,q)$ with $h_1\neq h_2$, we have
\[
\bil(c_{i,n,q} s_{i,n,q},g_{\q_j,m}^2) = (-1)^i (h_2-h_1) \delta_{n,m}\delta_{i,j}\qquad \text{for
  all $i,j\in \{1,2\}$ and $m,n\in\N$}
\]
for the normalized eigenfunctions $g_{\q_i,n}$ and $s_{i,n,q}$,
$c_{i,n,q}$ as defined above.
\end{lemma}
%
%\[
%|\gamma| = \frac{|h_1-h_2|}{\sqrt{(1+h_1^2)(1+h_2^2)}} 
%\] 

Using the alternative bilinear form 
\[
\begin{array}{rrcl}
\abil: &H^1([0,1],\R)× (\R^3×  L^2([0,1],\R))   & \longrightarrow & \R\\
   &   (f,(a,b,c,g)) &\mapsto &  -2\int_0^1f'g \di x +f(1)b +f(1)c+f(0)a
\end{array}\,,
\]
and integration by parts 
\[
\bil(f,g) = -2\int_0¹ f'g \di x +fg(1) -fg(0)\,,
\]
we get the corresponding statement for the gradients $\lambda_{\q_j,m}$
\[
\abil(c_{i,n,q} s_{i,n,q},\grad \lambda_{\q_j,m}) = \bil(c_{i,n,q} s_{i,n,q},g_{\q_j,m}^2)=(-1)^i (h_2-h_1) \delta_{n,m}\delta_{i,j}\,.
\]

Finally, since $\abil$ is bounded by
\[
|\abil(f,(a,b,c,g))| \leq 2\|f\|_{H^1}\|g\|_{L^2} + \sqrt{2}\|f\|_{H^1}(|a|+|b|+|c| )\leq 2\|f\|_{H^1} (\|g\|_{L^2} + |a| + |b| +|c|)
\]
and, in particular, continuous in the 2nd component, we immediately get
the linear independence theorem.

\begin{theorem}\label{linind}
With the notations of the above lemma, the set of gradients of the eigenvalues
\[
 \big\{ \grad \lambda_{\q_i,n} | (i,n)\in \{1,2\} × \N\big\} 
\]
is linearly independent in $\R^3×L^2$.
\end{theorem}
\begin{proof}
Suppose for some fixed $(i,n)$ we have
\[
\grad \lambda_{\q_i,n}=\sum_{k\in \N} a_k\grad \lambda_k
\]
in $\R^3×   L^2$,
where $a_k\in\rz$ and $\grad \lambda_k=\grad \lambda_{\q_{j_k},m_k}$ with $(j_k,m_k)\neq(i,n)$. But
this would imply
\[
(-1)^i(h_2-h_1) = \abil (c_{i,n,q}s_{i,n,q},\grad \lambda_{\q_i,n}) = 
\abil \left(c_{i,n,q}s_{i,n,q},\sum_{k\in\N} a_k\grad \lambda_k\right) = \sum_{k\in\N}\abil (c_{i,n,q}s_{i,n,q},a_k\grad \lambda_k) =0 
\] %
\end{proof}

%    Bibliographies can be prepared with BibTeX using amsplain,
%    amsalpha, or (for "historical" overviews) natbib style.
\bibliographystyle{amsplain}
\bibliography{inp}
%    Insert the bibliography data here.

\end{document}